\documentclass[reqno]{amsart}
\usepackage{amssymb,amsmath,enumerate,mathrsfs}

\numberwithin{equation}{section}

\newtheorem{theorem}[equation]{Theorem}
\newtheorem{proposition}[equation]{Proposition}
\newtheorem{lemma}[equation]{Lemma}
\newtheorem{corollary}[equation]{Corollary}

\theoremstyle{definition}
\newtheorem{definition}[equation]{Definition}
\newtheorem{remark}[equation]{Remark}

\def\C{\mathbb C}
\def\HT{(\mathcal {HT})}
\def\L{\mathscr L}
\def\N{\mathbb N}
\def\R{\mathbb R}
\def\S{\mathscr S}
\def\x{\mathbf x}
\def\V{\mathcal V}
\def\Mbar{\overline{M}}
\def\RR{\mathcal R}
\def\Dom{\mathcal D}

\DeclareMathOperator{\sym}{ \sigma\!\!\!\sigma}

\def\sc{\textup{sc}}
\def\scsym{\,{}^{\sc}\!\sym}
\def\scT{\,{}^{\sc} T}
\def\scpi{\,{}^{\sc}\hspace{-1.5pt}\pi}
\def\scH{{}^{\sc}H}
\def\scPsi{{}^{\sc}\Psi}
\def\cl{\textup{cl}}
\def\eps{\varepsilon}
\def\m{\mathfrak m}
\def\open#1{{\mathring{#1}}}
\def\st{;\;}
\DeclareMathOperator{\scDiff}{{}^{\sc}Diff}
\DeclareMathOperator{\SP}{SP}
\DeclareMathOperator{\spec}{spec}
\DeclareMathOperator{\op}{op}

\begin{document}

\nocite*

\title[$\RR$-boundedness, pseudodiff. operators, and maximal regularity]{$\RR$-boundedness, pseudodifferential operators, and maximal
regularity for some classes of partial differential operators}
\author{Robert Denk}
\address[Robert Denk]{Fachbereich Mathematik und Statistik \\ Fach D 193 \\ Universit{\"a}t Konstanz \\ D-78457 Konstanz \\ Germany}
\author{Thomas Krainer}
\address[Thomas Krainer]{Institut f{\"u}r Mathematik \\ Universit{\"a}t Potsdam \\ Am Neuen Palais~10 \\ D-14469 Potsdam \\ Germany}
\address[Thomas Krainer]{Department of Mathematics and Statistics \\ The Pennsylvania State
University \\ Altoona Campus \\ 3000 Ivyside Park \\ Altoona, PA 16601 \\ U.S.A.}
\begin{abstract}
It is shown that an elliptic scattering operator $A$ on a compact manifold with boundary with
operator valued coefficients in the morphisms of a bundle of Banach spaces of class
$\HT$ and Pisier's property $(\alpha)$ has maximal regularity (up to a spectral shift),
provided that the spectrum of the principal symbol of $A$ on the scattering cotangent bundle
avoids the right half-plane. This is accomplished by representing the resolvent
in terms of pseudodifferential operators with $\RR$-bounded symbols, yielding by an iteration
argument the $\RR$-boundedness of $\lambda(A-\lambda)^{-1}$ in $\Re(\lambda) \geq \gamma$ for some
$\gamma \in \R$. To this end, elements of a symbolic and operator calculus of pseudodifferential
operators with $\RR$-bounded symbols are introduced. The significance of this method
for proving maximal regularity results for partial differential operators is underscored
by considering also a more elementary situation of anisotropic elliptic operators on $\R^d$ with
operator valued coefficients.
\end{abstract}

\thanks{Version: July 26, 2006.}
\maketitle

\section{Introduction}

\noindent
A central question in the analysis of parabolic evolution equations is
whether a linear operator enjoys the property of maximal regularity.
An account on how maximal regularity can be applied to partial differential equations
is given in the survey paper by Pr{\"u}ss \cite{Pruess}, and a general
overview on developments in operator theory that are connected with maximal
regularity can be found in Kunstmann and Weis \cite{KunstmannWeis},
and the monograph \cite{DHP} by Denk, Hieber, and Pr{\"u}ss.

\begin{definition}\label{MRegDefIntro}
A closed and densely defined operator $A$ with domain $\Dom(A)$
in a Banach space $X$ is said to have \emph{maximal $L_p$--regularity}, if the associated
evolution equation
\begin{equation}\label{EvolGl}
\frac{d}{dt} - A : \open{W}^1_p([0,\infty),X)\cap L_p((0,\infty),\Dom(A)) \to L_p((0,\infty),X)
\end{equation}
is an isomorphism for some $1 < p < \infty$.
Here $\open{W}^1_p([0,\infty),X)$ consists of all $u \in L_p(\R,X)$ with
$u' \in L_p(\R,X)$, and $u$ supported in $[0,\infty)$.
\end{definition}

\noindent
We shall assume here and in the sequel that the Banach space $X$ if of class $\HT$,
i.e. the vector valued Hilbert transform is assumed to be continuous on $L_p(\R,X)$ (see also
Definition \ref{HTalpha}, more on such spaces can be found in \cite{AmannBook}, \cite{DHP}).
It is known that the condition of maximal regularity does not depend on $p$, i.e.
\eqref{EvolGl} is an isomorphism for all $1 < p < \infty$ once that this is the case for
some $p$.

If the Banach space $X$ is an $L_q$-space, an application of the closed graph
theorem to \eqref{EvolGl} reveals that the property of maximal regularity is connected
with proving optimal apriori $L_p$--$L_q$ estimates for solutions.
This gives a hint why maximal regularity is important in the theory of nonlinear partial
differential equations, because apriori estimates of
such kind for the linearized equation and a contraction principle readily imply local
existence of solutions to the nonlinear equation.
Recent work on optimal $L_p$--$L_q$ estimates for parabolic evolution equations that
rely on heat kernel estimates includes \cite{MazzucatoNistor}, see also \cite{HieberPruess}.

There are several approaches to prove maximal regularity for a given operator $A$.
One approach is to check whether $A$ admits a bounded $H^{\infty}$--calculus (or merely
bounded imaginary powers). A famous result by Dore and Venni \cite{DoreVenni} then
implies maximal regularity.
Another way involves only the resolvent of $A$ and relies on operator valued
Fourier multiplier theorems due to Weis \cite{Weis}:
$A$ has maximal $L_p$--regularity provided that $A - \lambda : \Dom(A) \to X$ is
invertible for all $\lambda \in \C$ with $\Re(\lambda) \geq 0$, and the resolvent
$\{\lambda(A-\lambda)^{-1}\st \Re(\lambda) \geq 0\} \subset \L(X)$ is $\RR$-bounded
(see Section \ref{sec-Rbdd} further below for the definition and properties of
$\RR$-bounded sets of operators).

For elliptic partial differential operators $A$, the way to get a hold on its powers
is via making use of the resolvent and a Seeley theorem, i.e.
a representation of $(A-\lambda)^{-1}$ in terms of pseudodifferential
operators (see the classical paper \cite{SeelCP}, or \cite{Shubin}). However, in the case of operators on noncompact and singular
manifolds, it may still
be quite difficult to actually pass from a corresponding Seeley theorem to the
powers (or $H^{\infty}$--calculus) of $A$.
The way to prove maximal regularity by establishing the $\RR$-boundedness
of resolvents seems to be more direct for differential operators.
The intuition is that a Seeley theorem alone should already be sufficient.

In this work, we follow the latter philosophy. We investigate pseudodifferential
operators depending on parameters that play the role of the spectral parameter, and
we prove that they give rise to families of continuous operators in Sobolev spaces
that are $\RR$-bounded with suitable bounds (see Theorems \ref{IterationRboundedness}
and \ref{IterationScattering}). Hence the representation of the resolvent
in terms of a pseudodifferential parametrix that depends on the spectral parameter
immediately yields the desired $\RR$-boundedness, thus maximal regularity.

To achieve this, we introduce and investigate in Sections \ref{sec-Rbdd} and \ref{sec-symbols}
some properties of a symbolic
calculus with operator valued symbols that satisfy symbol estimates in terms of $\RR$-bounds
rather than operator norm bounds, and we investigate associated
pseudodifferential operators. Surprisingly for us, not much seems to have been done
in this direction (see Strkalj \cite{Strkalj} for some results). One
observation in this context is that classical operator valued symbols that are modelled
on the operator norm (i.e. symbols that admit asymptotic expansions with homogeneous components) automatically satisfy the strong
$\RR$-bounded symbol estimates (Proposition \ref{classicalRbounded}). This is important
since classical symbols are the ones that appear naturally when constructing
parametrices for differential operators.
Theorem \ref{IterationRboundedness} on families of pseudodifferential operators can be interpreted as an iteration result: A pseudodifferential
operator depending on parameters with $\RR$-bounded symbol induces a family of continuous operators
in Sobolev spaces that satisfies suitable $\RR$-bounds. More precisely, this family
is itself an $\RR$-bounded operator valued symbol depending on the parameter.
The reader will certainly notice that much from Section \ref{sec-symbols} can (and should)
be generalized to wider classes of symbols and operators.

In Section \ref{Quasi} we illustrate the proposed method by proving the $\RR$-boundedness
of resolvents of parameter-dependent anisotropic elliptic operators $A$ on $\R^d$
in anisotropic Sobolev spaces (Theorem \ref{RbddResQuasi}).
A parameter-dependent parametrix of $A - \lambda$ can simply be constructed by the
standard method of symbolic inversion and a formal Neumann series argument, and the
results from Sections \ref{sec-Rbdd} and \ref{sec-symbols} give the desired conclusion.

Finally, in Section \ref{Scattering}, we consider the more advanced situation of elliptic
scattering operators $A$ on manifolds with boundary (see Melrose \cite{Melrose}). Assuming
the appropriate ellipticity condition on the principal scattering symbol of $A$ we show that
the resolvent is an $\RR$-bounded family in the scattering Sobolev spaces. This is again
achieved directly from a Seeley theorem, i.e. the resolvent $(A-\lambda)^{-1}$ is represented in terms of
a parameter-dependent parametrix in the scattering pseudodifferential calculus (Theorem \ref{IterationScattering}).

In both Sections \ref{Quasi} and \ref{Scattering} the differential operator $A$ in question is admitted
to have operator valued coefficients. The underlying Banach space (or bundle) is assumed to be of
class $\HT$ and to have Pisier's property $(\alpha)$ (see, e.g., \cite{CPSW}, \cite{Pisier}).
Besides of the observation that in this way the case of systems with infinitely many equations
and unknowns is included, operators of such kind have recently been investigated in coagulation and
fragmentation models (see Amann \cite{AmannCluster}, though, as pointed out in subsequent
work by Amann, in many cases one is interested in operator valued coefficients acting
on $L_1$-spaces that fail to be of class $\HT$).

\subsubsection*{Acknowledgement}

We would like to express our gratitude to Jan Pr{\"u}ss from the University of Halle for
several discussions, in particular for bringing to our attention the paper by Girardi and
Weis \cite{GirardiWeis}.

\section{Preliminaries on $\RR$-boundedness}\label{sec-Rbdd}

\begin{definition}\label{Rboundednessdef}
Let $X$ and $Y$ be Banach spaces. A subset ${\mathcal T} \subset \L(X,Y)$ is called
\emph{$\RR$-bounded}, if for some $1 \leq p < \infty$ and some constant $C_p \geq 0$
the inequality
\begin{equation}\label{Rboundedinequality}
\Bigl(\sum\limits_{\eps_1,\ldots,\eps_N \in \{-1,1\}}\Bigl\|\sum\limits_{j=1}^N\eps_j T_jx_j\Bigr\|^p\Bigr)^{1/p}
\leq C_p \Bigl(\sum\limits_{\eps_1,\ldots,\eps_N \in \{-1,1\}}\Bigl\|\sum\limits_{j=1}^N\eps_j x_j\Bigr\|^p\Bigr)^{1/p}
\end{equation}
holds for all choices of $T_1,\ldots,T_N \in {\mathcal T}$ and $x_1,\ldots,x_N \in X$, $N \in \N$.

The best constant
\begin{gather*}
C_p = \sup\Bigl\{\Bigl(\sum\limits_{\eps_1,\ldots,\eps_N \in \{-1,1\}}\Bigl\|\sum\limits_{j=1}^N\eps_j T_jx_j\Bigr\|^p\Bigr)^{1/p} \st
N \in \N,\: T_1,\ldots,T_N \in {\mathcal T}, \\
\Bigl(\sum\limits_{\eps_1,\ldots,\eps_N \in \{-1,1\}}\Bigl\|\sum\limits_{j=1}^N\eps_j x_j\Bigr\|^p\Bigr)^{1/p} = 1\Bigr\}
\end{gather*}
in \eqref{Rboundedinequality} is called the \emph{$\RR$-bound} of ${\mathcal T}$ and will be
denoted by $\RR({\mathcal T})$. By \emph{Kahane's inequality} \eqref{Kahaneinequality} the
notion of $\RR$-boundedness is independent of $1 \leq p < \infty$, and the $\RR$-bounds for
different values of $p$ are equivalent which is the justification for suppressing
$p$ from the notation.
\end{definition}

\begin{remark}
The following results related to $\RR$-bounded sets in $\L(X,Y)$ are well established in
the literature, see \cite{DHP, KunstmannWeis} and the references given there.
\begin{enumerate}[i)]
\item For Hilbert spaces $X$ and $Y$ the notion of $\RR$-boundedness reduces merely to boundedness.
\item \emph{Kahane's inequality}\/: For every Banach space $X$ and all values $1 \leq p,q < \infty$ there exist
constants $c,\: C > 0$ such that
\begin{equation}\label{Kahaneinequality}
\begin{gathered}
c\Bigl(2^{-N}\sum\limits_{\eps_1,\ldots,\eps_N \in \{-1,1\}}\Bigl\|\sum\limits_{j=1}^N\eps_j x_j\Bigr\|^q\Bigr)^{1/q}
\leq \Bigl(2^{-N}\sum\limits_{\eps_1,\ldots,\eps_N \in \{-1,1\}}\Bigl\|\sum\limits_{j=1}^N\eps_j x_j\Bigr\|^p\Bigr)^{1/p} \\
\leq C\Bigl(2^{-N}\sum\limits_{\eps_1,\ldots,\eps_N \in \{-1,1\}}\Bigl\|\sum\limits_{j=1}^N\eps_j x_j\Bigr\|^q\Bigr)^{1/q}
\end{gathered}
\end{equation}
for all choices $x_1,\ldots,x_N \in X$, $N \in \N$.
\item \emph{Kahane's contraction principle}\/: For all $\alpha_j,\: \beta_j \in \C$ with
$|\alpha_j| \leq |\beta_j|$, $j = 1,\ldots, N$, the inequality
\begin{equation}\label{Kahanecontraction}
\Bigl(\sum\limits_{\eps_1,\ldots,\eps_N \in \{-1,1\}}\Bigl\|\sum\limits_{j=1}^N\eps_j \alpha_jx_j\Bigr\|^p\Bigr)^{1/p}
\leq 2 \Bigl(\sum\limits_{\eps_1,\ldots,\eps_N \in \{-1,1\}}\Bigl\|\sum\limits_{j=1}^N\eps_j \beta_jx_j\Bigr\|^p\Bigr)^{1/p}
\end{equation}
holds for all choices $x_1,\ldots,x_N \in X$, $N \in \N$.

In particular, the set $\{\lambda I \st |\lambda| \leq R\} \subset \L(X)$ is $\RR$-bounded
for every $R > 0$.
\item For ${\mathcal T},\: {\mathcal S} \subset \L(X,Y)$ we have
\begin{equation}\label{Rboundsubadd}
\RR({\mathcal T}+{\mathcal S}) \leq \RR({\mathcal T}) + \RR({\mathcal S}).
\end{equation}
\item For Banach spaces $X,\: Y, \: Z$ and ${\mathcal T} \subset \L(Y,Z)$,
${\mathcal S} \subset \L(X,Y)$ we have
\begin{equation}\label{Rboundsubmult}
\RR({\mathcal T}{\mathcal S}) \leq \RR({\mathcal T})\RR({\mathcal S}).
\end{equation}
\item Let ${\mathcal T} \subset \L(X,Y)$, and let $\overline{\textup{aco}({\mathcal T})}$
be the closure of the absolute convex hull of ${\mathcal T}$ in the strong operator topology.
Then
\begin{equation}\label{Rboundacohull}
\RR\bigl(\overline{\textup{aco}({\mathcal T})}\bigr) \leq 2 \RR({\mathcal T}).
\end{equation}
\end{enumerate}
\end{remark}

\begin{definition}
Let $\Gamma$ be a set. Define $\ell_{\RR}^{\infty}(\Gamma,\L(X,Y))$ as the space
of all functions $f : \Gamma \to \L(X,Y)$ with $\RR$-bounded range and norm
\begin{equation}\label{Rboundnorm}
\|f\|:= \RR\bigl(f(\Gamma)\bigr).
\end{equation}
\end{definition}

\begin{proposition}\label{linftyRBanach}
$\bigl(\ell_{\RR}^{\infty}(\Gamma,\L(X,Y)),\|\cdot\|\bigr)$ is a Banach space. The
embedding
$$
\ell_{\RR}^{\infty}(\Gamma,\L(X,Y)) \hookrightarrow
\ell^{\infty}(\Gamma,\L(X,Y))
$$
into the Banach space $\ell^{\infty}(\Gamma,\L(X,Y))$ of all $\L(X,Y)$-valued functions
on $\Gamma$ with bounded range is a contraction.

The norm in $\ell^{\infty}_{\RR}$ is submultiplicative, i.e.
$$
\|f\cdot g\|_{\ell^{\infty}_{\RR}} \leq \|f\|_{\ell^{\infty}_{\RR}} \cdot \|g\|_{\ell^{\infty}_{\RR}}
$$
whenever the composition $f\cdot g$ makes sense, and we have $\|{\mathfrak 1}\|_{\ell^{\infty}_{\RR}} = 1$
for the constant map ${\mathfrak 1} \equiv \textup{Id}_X$.
\end{proposition}
\begin{proof}
Definiteness and homogeneity of the norm \eqref{Rboundnorm} are immediate consequences
of Definition \ref{Rboundednessdef}. The triangle inequality follows from \eqref{Rboundsubadd},
the submultiplicativity from \eqref{Rboundsubmult}.
Moreover, the embedding $\ell^{\infty}_{\RR} \hookrightarrow \ell^{\infty}$ is a contraction
because the $\RR$-bound of a set is always greater or equal to its operator norm bound.

It remains to show completeness. Let $(f_j)_j \subset \ell^{\infty}_{\RR}$ be a Cauchy sequence.
Thus $(f_j)_j$ is also a Cauchy sequence in $\ell^{\infty}(\Gamma,\L(X,Y))$, and there exists
$f \in \ell^{\infty}(\Gamma,\L(X,Y))$ with $\|f_j - f\|_{\ell^{\infty}} \to 0$ as $j \to \infty$.
Let $\eps > 0$, and let $N(\eps) \in \N$ be such that $\|f_j-f_k\|_{\ell^{\infty}_{\RR}} \leq \eps$
for $j,k \geq N(\eps)$. In view of Definition \ref{Rboundednessdef} this implies that
\begin{equation}\label{linfrcomplete}
\Bigl(\sum\limits_{\eps_1,\ldots,\eps_N \in \{-1,1\}}\Bigl\|\sum\limits_{i=1}^N\eps_i (f_j(\gamma_i)x_i-f_k(\gamma_i)x_i)\Bigr\|^p\Bigr)^{1/p} \leq \eps
\end{equation}
for all finite collections $x_1,\ldots,x_N \in X$ with
$\Bigl(\sum\limits_{\eps_1,\ldots,\eps_N \in \{-1,1\}}\Bigl\|\sum\limits_{j=1}^N\eps_i x_j\Bigr\|^p\Bigr)^{1/p} = 1$,
and all choices $\gamma_1,\ldots,\gamma_N \in \Gamma$.
Letting $k \to \infty$ in \eqref{linfrcomplete} gives
$$
\Bigl(\sum\limits_{\eps_1,\ldots,\eps_N \in \{-1,1\}}\Bigl\|\sum\limits_{i=1}^N\eps_i (f_j(\gamma_i)x_i-f(\gamma_i)x_i)\Bigr\|^p\Bigr)^{1/p} \leq \eps,
$$
and passing to the supremum over all possible choices implies
$\|f_j - f\|_{\ell^{\infty}_{\RR}} \leq \eps$ for $j \geq N(\eps)$. This shows
$f \in \ell^{\infty}_{\RR}$ and $f_j \to f$ with respect to $\|\cdot\|_{\ell^{\infty}_{\RR}}$,
and the proof is complete.
\end{proof}

\begin{proposition}\label{linftyproj}
We have
$$
\ell^{\infty}(\Gamma) \hat{\otimes}_{\pi} \L(X,Y) \subset \ell^{\infty}_{\RR}(\Gamma,\L(X,Y)).
$$
Recall that $\ell^{\infty}(\Gamma) \hat{\otimes}_{\pi} \L(X,Y)$ is realized as the space
of all functions $f: \Gamma \to \L(X,Y)$ that can be represented as
$$
f(\gamma) = \sum\limits_{j=1}^{\infty}\lambda_j f_j(\gamma)A_j
$$
with sequences $(\lambda_j)_j \in \ell^1(\N)$, and $f_j \to 0$ in $\ell^{\infty}(\Gamma)$
as well as $A_j \to 0$ in $\L(X,Y)$.
\end{proposition}
\begin{proof}
For $j \in \N$ the function $\lambda_jf_j \otimes A_j : \Gamma \to \L(X,Y)$ belongs to
$\ell^{\infty}_{\RR}(\Gamma,\L(X,Y))$ with norm
$$
\|\lambda_jf_j \otimes A_j\|_{\ell^{\infty}_{\RR}} \leq 2 \cdot |\lambda_j| \cdot \|f_j\|_{\ell^{\infty}(\Gamma)}\cdot \|A_j\|_{\L(X,Y)}
$$
in view of Kahane's contraction principle \eqref{Kahanecontraction}. Consequently,
$$
\sum\limits_{j=1}^{\infty}\|\lambda_jf_j \otimes A_j\|_{\ell^{\infty}_{\RR}} < \infty,
$$
and by completeness the series $f = \sum\limits_{j=1}^{\infty}\lambda_j f_j \otimes A_j$
converges in $\ell^{\infty}_{\RR}(\Gamma,\L(X,Y))$.
\end{proof}

As a consequence of Proposition \ref{linftyproj} we obtain the following corollary about
the $\RR$-boundedness of the range of certain very regular functions. Actually, much less
regularity is necessary to draw this conclusion, see, e.g., \cite{GirardiWeis}.
However, for our purposes the corollary is sufficient, and its proof is based on elementary
arguments.

\begin{corollary}\label{BilderRbounded}
\begin{enumerate}[i)]
\item Let $M$ be a smooth manifold, $K \subset M$ a compact subset, and $f \in C^{\infty}(M,\L(X,Y))$.
Then $f(K)$ is an $\RR$-bounded subset of $\L(X,Y)$.
\item Let $f \in \S(\R^n,\L(X,Y))$. Then the range $f(\R^n) \subset \L(X,Y)$ is $\RR$-bounded.
\end{enumerate}
\end{corollary}
\begin{proof}
The assertion follows from Proposition \ref{linftyproj} in view of
\begin{align*}
C^{\infty}(M,\L(X,Y)) &\cong C^{\infty}(M) \hat{\otimes}_{\pi} \L(X,Y), \\
\S(\R^n,\L(X,Y)) &\cong \S(\R^n) \hat{\otimes}_{\pi} \L(X,Y).
\end{align*}
\end{proof}

\section{Operator valued $\RR$-bounded symbols, and pseudodifferential operators on $\R^d$}\label{sec-symbols}

\noindent
In this section we consider special classes of anisotropic operator valued symbols
and associated pseudodifferential operators in $\R^d$ depending on parameters. In what
follows, let $n \in \N$ be total dimension of parameters and covariables.
Throughout this section we fix a vector $\vec{\ell} = (\ell_1,\ldots,\ell_n) \in \N^n$ of
positive integers which represents the anisotropy. For $\xi \in \R^n$ we denote
$$
|\xi|_{\vec{\ell}} = \Bigl(\sum\limits_{j=1}^n \xi_j^{2\pi_j}\Bigr)^{\frac{1}{2\ell_1\cdots\ell_n}}, \quad
\langle \xi \rangle_{\vec{\ell}} = \Bigl(1 + \sum\limits_{j=1}^n \xi_j^{2\pi_j}\Bigr)^{\frac{1}{2\ell_1\cdots\ell_n}}, \quad
\textup{where } \pi_j = \prod\limits_{i \neq j}\ell_i,
$$
and for a multi-index $\beta \in \N_0^n$ let $|\beta|_{\vec{\ell}} = \sum\limits_{j=1}^n\ell_j\beta_j$ be
its anisotropic length. We are aware that the notation $|\cdot|_{\vec{\ell}}$ is ambiguous,
on the other hand multi-indices and covectors (or parameters) are easily distinguishable by
the context where they appear.

Apparently, Peetre's inequality
$$
\langle \xi+\xi' \rangle_{\vec{\ell}}^s \leq 2^{|s|}\langle \xi \rangle_{\vec{\ell}}^s
\cdot\langle \xi' \rangle_{\vec{\ell}}^{|s|}
$$
holds for all $s \in \R$, and there exist constants $c,C > 0$ depending only on $\vec{\ell}$ and
$n$ such that
$$
c \langle \xi \rangle^{1/{\sum\limits_{j=1}^n\ell_j}} \leq \langle \xi \rangle_{\vec{\ell}} \leq
C \langle \xi \rangle^{\sum\limits_{j=1}^n 1/{\ell_j}},
$$
where as usual $\langle \xi \rangle = \bigl(1 + |\xi|^2\bigr)^{1/2}$ is the standard
regularized distance function.

\medskip

Let $X$ and $Y$ be Banach spaces, $\mu \in \R$, and let $S^{\mu;\vec{\ell}}(\R^n;X,Y)$ denote
the standard space of anisotropic $\L(X,Y)$-valued symbols on $\R^n$, i.e. the space of
all $a \in C^{\infty}(\R^n,\L(X,Y))$ such that
$$
\sup\limits_{\xi \in \R^n}\|\langle \xi \rangle_{\vec{\ell}}^{-\mu+|\beta|_{\vec{\ell}}}\partial_{\xi}^{\beta}a(\xi)\|_{\L(X,Y)} < \infty
$$
for all $\beta \in \N_0^n$.

We shall be mainly concerned with the following more restrictive symbol class:

\begin{definition}\label{symboldef}
A function $a \in C^{\infty}(\R^n,\L(X,Y))$ belongs to
$S^{\mu;\vec{\ell}}_{\RR}(\R^n;X,Y)$ if and only if for all $\beta \in \N_0^n$
\begin{equation}\label{symbolestimate}
|a|^{(\mu;\vec{\ell})}_{\beta}:= {\mathcal R}\bigl(\{\langle \xi \rangle_{\vec{\ell}}^{-\mu+|\beta|_{\vec{\ell}}}\partial_{\xi}^{\beta}a(\xi)\st \xi \in \R^n\}\bigr) < \infty.
\end{equation}
By Proposition \ref{linftyRBanach} and the usual arguments we obtain that
$S^{\mu;\vec{\ell}}_{{\RR}}(\R^n;X,Y)$ is a Fr{\'e}chet space in the topology
generated by the seminorms $|\cdot|^{(\mu;\vec{\ell})}_{\beta}$.
\end{definition}

\begin{lemma}\label{symbolelementary}
\begin{enumerate}[i)]
\item For $\beta \in \N_0^n$ differentiation
$$
\partial_{\xi}^{\beta} :
S^{\mu;\vec{\ell}}_{\RR}(\R^n;X,Y) \to S^{\mu-|\beta|_{\vec{\ell}};\vec{\ell}}_{\RR}(\R^n;X,Y)
$$
is continuous.
\item The embedding
$$
S^{\mu;\vec{\ell}}_{{\RR}}(\R^n;X,Y) \hookrightarrow S^{\mu';\vec{\ell}}_{{\RR}}(\R^n;X,Y)
$$
is continuous for $\mu \leq \mu'$.
\item For Banach spaces $X$, $Y$, and $Z$ the multiplication (pointwise composition)
$$
S^{\mu;\vec{\ell}}_{{\RR}}(\R^n;Y,Z)\times S^{\mu';\vec{\ell}}_{{\RR}}(\R^n;X,Y) \to
S^{\mu+\mu';\vec{\ell}}_{{\RR}}(\R^n;X,Z)
$$
is bilinear and continuous.
\item The embedding
$$
S^{\mu;\vec{\ell}}_{{\RR}}(\R^n;X,Y) \hookrightarrow S^{\mu;\vec{\ell}}(\R^n;X,Y)
$$
into the standard class of anisotropic operator valued symbols is continuous.
\item The space of scalar symbols $S^{\mu;\vec{\ell}}(\R^n)$ embeds into
$S^{\mu;\vec{\ell}}_{{\RR}}(\R^n;X,X)$ via $a(\xi) \mapsto a(\xi)\textup{Id}_X$.
\end{enumerate}
\end{lemma}
\begin{proof}
i) and iv) are evident, ii) and v) follow from Kahane's contraction principle \eqref{Kahanecontraction},
and iii) is straightforward from the Leibniz rule and Proposition \ref{linftyRBanach}
(submultiplicativity and subadditivity of the norm $\|\cdot\|_{\ell^{\infty}_{\RR}}$).
\end{proof}

\begin{lemma}\label{Smoothingcoincide}
We have $S^{-\infty}_{\RR}(\R^n;X,Y) = S^{-\infty}(\R^n;X,Y)$, i.e.
\begin{equation}\label{Smoothcoincide}
\bigcap\limits_{\mu \in \R}S^{\mu;\vec{\ell}}_{\RR}(\R^n;X,Y) = \bigcap\limits_{\mu\in\R}S^{\mu;\vec{\ell}}(\R^n;X,Y),
\end{equation}
and this space does not depend on the anisotropy $\vec{\ell}$.
\end{lemma}
\begin{proof}
By Lemma \ref{symbolelementary} we have $S^{-\infty}_{\RR} \subset S^{-\infty}$. On the other hand, we may write
$$
S^{-\infty}(\R^n;X,Y) \cong \S(\R^n)\hat{\otimes}_{\pi}\L(X,Y),
$$
and thus we obtain $S^{-\infty} \subset S_{\RR}^{-\infty}$ from Proposition \ref{linftyproj}.
\end{proof}

\begin{definition}\label{asympdef}
Let $a_j \in S^{\mu_j;\vec{\ell}}_{{\RR}}(\R^n;X,Y)$ with $\mu_j \to -\infty$, and let
$\overline{\mu} = \max\mu_j$. For a symbol $a \in S^{\overline{\mu};\vec{\ell}}_{{\RR}}(\R^n;X,Y)$
write $a \underset{\RR}{\sim} \sum\limits_{j=1}^{\infty}a_j$ if for all $K \in \R$ there exists $N(K)$ such that
$$
a - \sum\limits_{j=1}^{N}a_j \in S^{K;\vec{\ell}}_{{\RR}}(\R^n;X,Y)
$$
for $N > N(K)$.

Recall that the standard notion of asymptotic expansion $a \sim \sum\limits_{j=1}^{\infty}a_j$
means that for all $K \in \R$ there exists $N(K)$ such that
$$
a - \sum\limits_{j=1}^{N}a_j \in S^{K;\vec{\ell}}(\R^n;X,Y)
$$
for $N > N(K)$.
\end{definition}

\begin{proposition}\label{asympexist}
Let $a_j \in S^{\mu_j;\vec{\ell}}_{{\RR}}(\R^n;X,Y)$ with $\mu_j \to -\infty$, and let
$\overline{\mu} = \max\mu_j$. Then there exists $a \in S^{\overline{\mu};\vec{\ell}}_{{\RR}}(\R^n;X,Y)$
with $a \underset{\RR}{\sim} \sum\limits_{j=1}^{\infty}a_j$.
\end{proposition}
\begin{proof}
The proof is based on the usual Borel argument. Let $\chi \in C^{\infty}(\R^n)$ be a
function with $\chi \equiv 0$ near the origin and $\chi \equiv 1$ near infinity, and
define
$$
\chi_{\theta}(\xi):= \chi\Bigl(\frac{\xi_1}{\theta^{\ell_1}},\ldots,\frac{\xi_n}{\theta^{\ell_n}}\Bigr)
$$
for $\theta \geq 1$.
Then the family $\{\chi_{\theta}\st \theta \geq 1\} \subset S^{0;\vec{\ell}}(\R^n)$ is bounded.
To see this assume that $\chi \equiv 1$ for $|\xi|_{\vec{\ell}} \geq c$ and
$\chi \equiv 0$ for $|\xi|_{\vec{\ell}} \leq \frac{1}{c}$. Then, for $\beta \neq 0$,
$$
\partial_{\xi}^{\beta}\chi_{\theta}(\xi) = \theta^{-|\beta|_{\vec{\ell}}}
\bigl(\partial_{\xi}^{\beta}\chi\bigr)\Bigl(\frac{\xi_1}{\theta^{\ell_1}},\ldots,\frac{\xi_n}{\theta^{\ell_n}}\Bigr) \neq 0
$$
at most for $\frac{1}{c}|\xi|_{\vec{\ell}} \leq \theta \leq c|\xi|_{\vec{\ell}}$.

Let $p \in S^{\mu;\vec{\ell}}_{{\RR}}(\R^n;X,Y)$, $\mu \in \R$, and define
$p_{\theta}(\xi) = \chi_{\theta}(\xi)p(\xi)$, $\theta \geq 1$.
By Lemma \ref{symbolelementary} the set $\{p_{\theta}\st \theta \geq 1\} \subset S^{\mu;\vec{\ell}}_{{\RR}}(\R^n;X,Y)$
is bounded. Let $\mu' > \mu$ and $\eps > 0$. Fix $R(\eps) \geq 1$ such that for $\theta \geq R(\eps)$ we have
$\chi_{\theta}(\xi) \neq 0$ at most for $\xi \in \R^n$ with $\langle \xi \rangle_{\vec{\ell}}^{\mu'-\mu} \geq \frac{2}{\eps}$.
Consequently, for such $\theta$ we obtain from Kahane's contraction principle \eqref{Kahanecontraction}
\begin{gather*}
\Bigl(\sum\limits_{\eps_1,\ldots,\eps_N \in \{-1,1\}}\Bigl\|\sum\limits_{j=1}^N\eps_j 
\langle \xi_j \rangle_{\vec{\ell}}^{-\mu'+|\beta|_{\vec{\ell}}}\partial_{\xi}^{\beta}p_{\theta}(\xi_j)x_j\Bigr\|^p\Bigr)^{1/p} \\
\leq \eps
\Bigl(\sum\limits_{\eps_1,\ldots,\eps_N \in \{-1,1\}}\Bigl\|\sum\limits_{j=1}^N\eps_j 
\langle \xi_j \rangle_{\vec{\ell}}^{-\mu+|\beta|_{\vec{\ell}}}\partial_{\xi}^{\beta}p_{\theta}(\xi_j)x_j\Bigr\|^p\Bigr)^{1/p} \\
\leq \eps |p_{\theta}|^{(\mu;\vec{\ell})}_{\beta}
\Bigl(\sum\limits_{\eps_1,\ldots,\eps_N \in \{-1,1\}}\Bigl\|\sum\limits_{j=1}^N\eps_j 
x_j\Bigr\|^p\Bigr)^{1/p}
\end{gather*}
which shows that $|p_{\theta}|^{(\mu';\vec{\ell})}_{\beta} \leq \eps |p_{\theta}|^{(\mu;\vec{\ell})}_{\beta}$
for $\theta \geq R(\eps)$. The boundedness of $\{p_{\theta}\st \theta \geq 1\}$ in
$S^{\mu;\vec{\ell}}_{\RR}$ now implies that $p_{\theta} \to 0$ in $S^{\mu';\vec{\ell}}_{\RR}$ as $\theta \to \infty$.

Now we proceed with the construction of the symbol $a$. We may assume without loss of generality
that $\mu_j > \mu_{j+1}$ for $j \in {\mathbb N}$. For each $j \in {\mathbb N}$
let $p^j_1 \leq p^j_2 \leq \ldots$ be an increasing fundamental system of seminorms for the
topology of $S^{\mu_j;\vec{\ell}}_{\RR}$.
Pick a sequence $1 \leq c_k^1 < c_{k+1}^1 \to \infty$ such that
$p_k^1(a_{k,\theta}) < 2^{-k}$ for $k > 1$ and all $\theta \geq c_k^1$. We proceed by induction
and construct successively subsequences $(c_k^{j})_k$ of $(c_k^{j-1})_k$ such that
$p_k^{j}(a_{k,\theta}) < 2^{-k}$ for $k > j$ and all $\theta \geq c_k^{j}$. Let
$c_k:= c_k^k$ be the diagonal sequence. Then $p_k^{j}(a_{k,c_k}) < 2^{-k}$ for $k > j$
which shows that the series $\sum\limits_{k=j}^{\infty}a_{k,c_k}$ is unconditionally convergent
in $S_{\RR}^{\mu_j;\vec{\ell}}$. Now let $a:= \sum\limits_{k=1}^{\infty}a_{k,c_k}$. Then
$$
a - \sum\limits_{k=1}^Na_k = \sum\limits_{k=N+1}^{\infty}a_{k,c_k} + \sum\limits_{k=1}^N(a_{k,c_k}-a_k).
$$
As every summand $a_{k,c_k} - a_k$ is compactly supported in $\xi$, we obtain from
Lemma \ref{Smoothingcoincide} that
$$
\sum\limits_{k=1}^N(a_{k,c_k}-a_k) \in S^{-\infty}_{{\RR}}(\R^n;X,Y),
$$
and consequently $a - \sum\limits_{k=1}^Na_k \in S^{\mu_{N+1};\vec{\ell}}_{\RR}$ as desired.
\end{proof}

The combination of the possibility to carry out asymptotic expansions within the classes $S^{\mu;\vec{\ell}}_{\RR}$ and
the identity $S^{-\infty} = S^{-\infty}_{\RR}$ are very useful for proving that the $\RR$-bounded symbol
classes remain preserved under manipulations of the symbolic calculus.
Let us formulate this more precisely:

\begin{lemma}\label{asymptoticprinciple}
Let $\overline{\mu} = \max\limits_{j=1}^{\infty}\mu_j$, where $\mu_j \to -\infty$ as $j \to \infty$, and
let $a \in S^{\overline{\mu};\vec{\ell}}$ with $a \sim \sum\limits_{j=1}^{\infty}a_j$. Suppose that the summands
in the asymptotic expansion satisfy $a_j \in S^{\mu_j;\vec{\ell}}_{\RR}$ for every $j \in \N$. Then
$a \in S^{\overline{\mu};\vec{\ell}}_{\RR}$, and $a \underset{\RR}{\sim} \sum\limits_{j=1}^{\infty}a_j$.
\end{lemma}
\begin{proof}
According to Proposition \ref{asympexist} there exists a symbol $c \in S^{\overline{\mu};\vec{\ell}}_{\RR}$ with
$c \underset{\RR}{\sim} \sum\limits_{j=1}^{\infty}a_j$. In particular we have
$c \sim \sum\limits_{j=1}^{\infty}a_j$, and consequently $a - c \in S^{-\infty} = S^{-\infty}_{\RR}$
by Lemma \ref{Smoothingcoincide}. Thus $a \in S^{\overline{\mu};\vec{\ell}}_{\RR}$ with
$a \underset{\RR}{\sim} c \underset{\RR}{\sim} \sum\limits_{j=1}^{\infty}a_j$ as desired.
\end{proof}

We are going to make use of Lemma \ref{asymptoticprinciple} to show that the standard classes
of classical anisotropic symbols are a subspace of the $\RR$-bounded ones, see Proposition \ref{classicalRbounded}
further below. To this end, recall the definition of the spaces of classical symbols:

\begin{definition}
\begin{enumerate}[i)]
\item A function $a \in C^{\infty}\bigl(\R^n\setminus\{0\},\L(X,Y)\bigr)$ is called \emph{(anisotropic)
homogeneous of degree $\mu \in \R$} iff
$$
a\bigl(\varrho^{\ell_1}\xi_1,\ldots,\varrho^{\ell_n}\xi_n\bigr) = \varrho^{\mu}a(\xi)
$$
for $\varrho > 0$ and $\xi \in \R^n\setminus\{0\}$.
\item The space $S_{\cl}^{\mu;\vec{\ell}}(\R^n;X,Y)$ of \emph{classical symbols} of order $\mu$
consists of all $a \in S^{\mu;\vec{\ell}}(\R^n;X,Y)$ such that there exists a sequence
$a_{(\mu-j)}$ of anisotropic homogeneous functions of degree $\mu - j$, $j \in \N_0$, such
that for some excision function $\chi \in C^{\infty}(\R^n)$ with $\chi \equiv 0$ near the origin
and $\chi \equiv 1$ near infinity $a(\xi) \sim \sum\limits_{j=0}^{\infty}\chi(\xi)a_{(\mu-j)}(\xi)$.
\end{enumerate}
Recall that $S_{\cl}^{\mu;\vec{\ell}}(\R^n;X,Y)$ is a Fr{\'e}chet space in the projective
topology with respect to the mappings
\begin{align*}
a &\mapsto a(\xi) - \sum\limits_{j=0}^N\chi(\xi)a_{(\mu-j)}(\xi) \in S^{\mu-N-1;\vec{\ell}}(\R^n;X,Y), \qquad N=-1,0,1,\ldots, \\
a &\mapsto a_{(\mu-j)}|_{\{|\xi|_{\vec{\ell}}=1\}} \in C^{\infty}\bigl(\{|\xi|_{\vec{\ell}}=1\},\L(X,Y)\bigr), \qquad j=0,1,2,\ldots.
\end{align*}
Note that the function $\xi \mapsto \langle \xi \rangle_{\vec{\ell}}^{\mu}$ is a classical
scalar symbol of order $\mu \in \R$.
\end{definition}

\begin{proposition}\label{classicalRbounded}
$S^{\mu;\vec{\ell}}_{\cl}(\R^n;X,Y) \hookrightarrow S^{\mu;\vec{\ell}}_{\RR}(\R^n;X,Y)$.
\end{proposition}
\begin{proof}
Let $a \in S^{\mu;\vec{\ell}}_{\cl}(\R^n;X,Y)$. In view of
$a(\xi) \sim \sum\limits_{j=0}^{\infty}\chi(\xi)a_{(\mu-j)}(\xi)$ and Lemma \ref{asymptoticprinciple}
it suffices to show that each summand $\chi(\xi)a_{(\mu-j)}(\xi)$ in this asymptotic expansion
belongs to $S^{\mu-j;\vec{\ell}}_{\RR}(\R^n;X,Y)$.

Consider therefore the function $p(\xi) = \chi(\xi)a_{(\mu)}(\xi)$.
Note that $\partial_{\xi}^{\beta}a_{(\mu)}(\xi)$ is anisotropic homogeneous of degree
$\mu-|\beta|_{\vec{\ell}}$. Hence, by the Leibniz rule and the fact that compactly supported smooth
functions in $\xi$ belong to $S^{-\infty}_{\RR}$ by Lemma \ref{Smoothingcoincide},
we conclude that it is sufficient to prove that the seminorm $|p|^{(\mu;\vec{\ell})}_{0} < \infty$,
see \eqref{symbolestimate}.
We have
$$
\langle \xi \rangle_{\vec{\ell}}^{-\mu}\bigl(\chi(\xi)a_{(\mu)}(\xi)\bigr)
= \chi(\xi)a\Bigl(\frac{\xi_1}{\langle \xi \rangle_{\vec{\ell}}^{\ell_1}},\ldots,
\frac{\xi_n}{\langle \xi \rangle_{\vec{\ell}}^{\ell_n}}\Bigr),
$$
and thus this operator family is $\RR$-bounded for $\xi \in \R^n$ in view of
Corollary \ref{BilderRbounded}.
\end{proof}

\begin{remark}
Let $\emptyset \neq \Lambda \subset \R^n$. The symbol space $S_{\RR}^{\mu;\vec{\ell}}(\Lambda;X,Y)$
is defined as the space of restrictions of the $\RR$-bounded symbol class on $\R^n$
to $\Lambda$ endowed with the quotient topology, and analogously we define the classes of
(classical) ordinary symbols on $\Lambda$. In this way, the results of this section carry
over immediately to symbols on $\Lambda$.
\end{remark}

\subsection*{Pseudodifferential operators}

Let $n = d + q$, and we change the notation slightly so as to consider covariables
$\xi \in \R^d$ and parameters $\lambda \in \Lambda \subset \R^q$.
Let $\vec{\ell} = (\vec{\ell}',\vec{\ell}'') \in \N^{d+q}$ be the vector that determines
the anisotropy of covariables and parameters.

With a symbol $a(x,\xi,\lambda) \in S^0_{\cl}(\R_x^d,S_{\cl}^{\mu;\vec{\ell}}(\R_{\xi}^d\times\Lambda;X,Y))$
we associate a family of pseudodifferential operators $\op_x(a)(\lambda) : \S(\R^d,X) \to \S(\R^d,Y)$ as
usual via
$$
\op_x(a)(\lambda)u(x) = (2\pi)^{-d}\int\limits_{\R^d}e^{ix\xi}a(x,\xi,\lambda)\hat{u}(\xi)\,d\xi,
$$
where $\hat{u}(\xi) = \int\limits_{\R^d}e^{-iy\xi}u(y)\,dy$ is the Fourier transform
of the function $u \in \S(\R^d,X)$.
Here and in what follows, much less requirements on the behaviour at infinity
of the $x$-dependence of the symbol $a(x,\xi,\lambda)$ are necessary. However, for our
purposes this (rather strong) condition is sufficient, and the arguments become considerably
simpler. We first observe the following

\begin{lemma}\label{globalautomaticallyRbounded}
Let $a(x,\xi,\lambda) \in S^0_{\cl}(\R_x^d,S_{\cl}^{\mu;\vec{\ell}}(\R_{\xi}^d\times\Lambda;X,Y))$.
Then for every $\alpha,\beta$ the set
\begin{equation}\label{alphabetaset}
\bigl\{\langle x \rangle^{|\alpha|}\langle \xi,\lambda \rangle_{\vec{\ell}}^{-\mu+|\beta|_{\vec{\ell}}}\partial_x^{\alpha}\partial_{(\xi,\lambda)}^{\beta}a(x,\xi,\lambda) \st x \in \R^d,\,(\xi,\lambda) \in \R^d\times\Lambda\bigr\}
\end{equation}
is $\RR$-bounded in $\L(X,Y)$.
\end{lemma}
\begin{proof}
We have
$$
S^0_{\cl}(\R_x^d,S_{\cl}^{\mu;\vec{\ell}}(\R_{\xi}^d\times\Lambda;X,Y)) \cong
S^0_{\cl}(\R_x^d) \hat{\otimes}_{\pi} S_{\cl}^{\mu;\vec{\ell}}(\R_{\xi}^d\times\Lambda;X,Y),
$$
and thus we may write $a(x,\xi,\lambda) = \sum\limits_{j=1}^{\infty}\lambda_j f_j(x)a_j(\xi,\lambda)$
with $(\lambda_j)_j \in \ell^1(\N)$ and $f_j \to 0$ in $S^{0}_{\cl}$ and
$a_j \to 0$ in $S^{\mu;\vec{\ell}}_{\cl}$.

Let $|\cdot|_{\alpha,\beta}$ be the seminorm of a symbol defined by the $\RR$-bound
of the set \eqref{alphabetaset} (for any given function $a(x,\xi,\lambda)$). By Kahane's
contraction principle \eqref{Kahanecontraction} and Proposition \ref{classicalRbounded} we obtain
that for every $\alpha,\beta$ there exists a constant $C > 0$ and
continuous seminorms $q$ on $S^{0}_{\cl}$ and $p$ on $S^{\mu;\vec{\ell}}_{\cl}$ such
that $|\lambda_jf_j(x)a_j(\xi,\lambda)|_{\alpha,\beta} \leq C|\lambda_j|q(f_j)p(a_j)$ for
all $j \in \N$. Consequently, the series
$$
\sum\limits_{j=1}^{\infty}|\lambda_jf_j(x)a_j(\xi)|_{\alpha,\beta} < \infty
$$
for all $\alpha,\beta$ which implies the assertion in view of Proposition \ref{linftyRBanach}.
\end{proof}

It follows from standard results in the theory of pseudodifferential operators that the class
of operators we are considering is invariant under composition on rapidly decreasing functions, i.e. if $X$,
$Y$, and $Z$ are Banach spaces and
\begin{align*}
a(x,\xi,\lambda) &\in S^0_{\cl}(\R_x^d,S_{\cl}^{\mu_1;\vec{\ell}}(\R_{\xi}^d\times\Lambda;Y,Z)), \\
b(x,\xi,\lambda) &\in S^0_{\cl}(\R_x^d,S_{\cl}^{\mu_2;\vec{\ell}}(\R_{\xi}^d\times\Lambda;X,Y)),
\end{align*}
then the composition $\op_x(a)(\lambda)\circ\op_x(b)(\lambda) : \S(\R^d,X) \to \S(\R^d,Z)$
is again a pseudodifferential operator $\op_x(a{\#}b)(\lambda)$ with symbol
$$
(a{\#}b)(x,\xi,\lambda) \in S^0_{\cl}(\R_x^d,S_{\cl}^{\mu_1+\mu_2;\vec{\ell}}(\R_{\xi}^d\times\Lambda;X,Z)).
$$
The mapping $(a,b) \mapsto a{\#}b$ is bilinear and continuous in the symbol topology, and
$a{\#}b \sim \sum\limits_{\alpha\in\N_0^n}\frac{1}{\alpha!}\bigl(\partial_{\xi}^{\alpha}a\bigr)\bigl(D_{x}^{\alpha}b\bigr)$
in the sense that for every $K \in \R$ there is an $N(K)$ such that
$$
a{\#}b - \sum\limits_{|\alpha| \leq N}\frac{1}{\alpha!}\bigl(\partial_{\xi}^{\alpha}a\bigr)\bigl(D_{x}^{\alpha}b\bigr)
\in S^{K}(\R_x^d,S^{K;\vec{\ell}}(\R_{\xi}^d\times\Lambda;X,Z))
$$
for all $N > N(K)$.

\medskip

We are mainly interested in the mapping properties and dependence on the parameter
$\lambda \in \Lambda$ of the operators $\op_x(a)(\lambda)$ in (anisotropic) vector valued
Bessel potential spaces. To this end, recall the following

\begin{definition}\label{HsDef}
For $s \in \R$, $1 < p < \infty$, let $H_p^{s;\vec{\ell}'}(\R^d,X)$ be the completion of $\S(\R^d,X)$
with respect to the norm
$$
\|u\|_{H_p^{s;\vec{\ell}'}}:= \|\op\bigl(\langle \xi \rangle_{\vec{\ell}'}^s\textup{Id}_X\bigr)u\|_{L_p(\R^d,X)}.
$$
\end{definition}

In order to proceed further we have to impose conditions on the Banach spaces involved
(\cite{AmannBook, CPSW, DHP, GirardiWeis, Pisier}):

\begin{definition}\label{HTalpha}
A Banach space $X$
\begin{enumerate}[i)]
\item is of \emph{class $\HT$} if the Hilbert transform is continuous in $L_p(\R,X) \to L_p(\R,X)$ for
some (all) $1 < p < \infty$. Recall that the Hilbert transform is the Fourier multiplier with symbol
$\chi_{[0,\infty)}\textup{Id}_X$, where $\chi_{[0,\infty)}$ is the characteristic function
on the interval $[0,\infty)$.

Banach spaces of class $\HT$ are often called UMD-spaces, and it is worth mentioning
that they are necessarily reflexive.
\item has Pisier's \emph{property $(\alpha)$} if for some $1 \leq p < \infty$ there exists a
constant $C_p > 0$ such that the inequality
\begin{gather*}
\Bigl(\sum\limits_{\eps_1,\ldots,\eps_N \in \{-1,1\}}\sum\limits_{\eps'_1,\ldots,\eps'_N \in \{-1,1\}}
\Bigl\|\sum\limits_{j,k=1}^{N}\eps_{j}\eps'_{k}\alpha_{jk}x_{jk}\Bigr\|^p\Bigr)^{1/p} \\
\leq C_p
\Bigl(\sum\limits_{\eps_1,\ldots,\eps_N \in \{-1,1\}}\sum\limits_{\eps'_1,\ldots,\eps'_N \in \{-1,1\}}
\Bigl\|\sum\limits_{j,k=1}^{N}\eps_{j}\eps'_{k}x_{jk}\Bigr\|^p\Bigr)^{1/p}
\end{gather*}
holds for all choices $x_{jk} \in X$ and $\alpha_{jk} \in \{-1,1\}$, $j,k = 1,\ldots,N$,
$N \in \N$.
\end{enumerate}
\end{definition}
Both $\HT$ and $(\alpha)$ are purely topological properties of a Banach space and
depend only on its isomorphism class. Some remarks about the permanence properties of these
conditions are in order:
\begin{enumerate}[a)]
\item Hilbert spaces are of class $\HT$ and have property $(\alpha)$.
\item (Finite) direct sums, closed subspaces, and quotients by closed subspaces of Banach
spaces of class $\HT$ are again of class $\HT$.
\item If $X$ is of class $\HT$ or has property $(\alpha)$, then so does
$L_p(\Omega,X)$, $1 < p < \infty$, for any $\sigma$-finite measure space $\Omega$.
\end{enumerate}
Granted this, we immediately obtain a wealth of examples for Banach spaces of class $\HT$
and property $(\alpha)$, most notably spaces of functions.
We will generally assume henceforth that all Banach spaces are of class $\HT$ and have property $(\alpha)$.

\medskip

The following specialization of a Fourier multiplier theorem by Girardi and Weis
is essential for us:

\begin{theorem}[\cite{GirardiWeis}, Theorem 3.2]\label{iteratedFouriermultiplier}
Let $X$ and $Y$ be Banach spaces of class $\HT$ with property $(\alpha)$. Let
${\mathcal T} \subset \L(X,Y)$ be an $\RR$-bounded subset, and let
$$
{\mathcal M}({\mathcal T}):= \{m : \R^d\setminus\{0\} \to \L(X,Y) \st \xi^{\beta}\partial_{\xi}^{\beta}m(\xi) \in {\mathcal T}
\textup{ for } \xi \neq 0,\, \beta \leq (1,\ldots,1)\}.
$$
For $1 < p < \infty$ the set of associated Fourier multipliers
$$
\{{\mathcal F}_{\xi \to x}^{-1}m(\xi){\mathcal F}_{x \to \xi} \st m \in {\mathcal M}({\mathcal T})\} \subset
\L\bigl(L_p(\R^d,X),L_p(\R^d,Y)\bigr)
$$
is then well defined \emph{(}i.e. the Fourier multipliers are continuous in $L_p(\R^d,X) \to L_p(\R^d,Y)$\emph{)} and $\RR$-bounded
by $C\cdot\RR({\mathcal T})$ with a constant $C \geq 0$ depending only on $X$, $Y$, $p$, and the dimension $d$.
\end{theorem}

\begin{remark}
From Theorem \ref{iteratedFouriermultiplier} we get in particular that the embedding
$$
H_p^{s;\vec{\ell}'}(\R^d,X) \hookrightarrow H_p^{t;\vec{\ell}'}(\R^d,X)
$$
is well defined and continuous for $s \geq t$ and $1 < p < \infty$ for any Banach space $X$ of
class $\HT$ satisfying property $(\alpha)$.

Moreover, for values $s \in \N$ such that $\ell'_j\; |\; s$ for all $j=1,\ldots,d$,
where $\vec{\ell}' = (\ell'_1,\ldots,\ell'_d) \in \N^d$,
the space $H_p^{s;\vec{\ell}'}(\R^d,X)$ coincides with the anisotropic Sobolev space
$$
W_p^{s;\vec{\ell}'}(\R^d,X) = \{u \in \S'(\R^d,X)\st \partial_x^{\alpha}u \in L_p(\R^d,X)
\textup{ for all } |\alpha|_{\vec{\ell}'} \leq s\}.
$$
The latter follows by the usual argument: First, the Fourier multipliers ${\mathcal F}^{-1}_{\xi \to x}\xi^{\alpha}\cdot\textup{Id}_X{\mathcal F}_{x \to \xi}$
are continuous in $H_p^{s;\vec{\ell}'}(\R^d,X) \to L_p(\R^d,X)$ for all $|\alpha|_{\vec{\ell}'} \leq s$,
which shows the inclusion $H_p^{s;\vec{\ell}'}(\R^d,X) \subset W_p^{s;\vec{\ell}'}(\R^d,X)$.
On the other hand, there exist classical scalar symbols $\varphi_j \in S^{0;\vec{\ell}'}_{\cl}(\R^d_{\xi})$,
$j = 1,\ldots,d$, such that the function
$$
m(\xi) = 1 + \sum\limits_{j=1}^d\varphi_j(\xi)\xi_j^{s/{\ell_j'}} \in S^{s;\vec{\ell}'}_{\cl}(\R^d_{\xi})
$$
satisfies $m(\xi) \geq c\langle \xi \rangle_{\vec{\ell}'}^s$ for all $\xi \in \R^d$ with some
constant $c > 0$, e.g. choose
$\varphi_j(\xi) = \chi(\xi)|\xi|_{\vec{\ell'}}^{-s}\xi_j^{s/{\ell_j'}}$, where $\chi \in C^{\infty}(\R^d_{\xi})$
is a suitable excision function of the origin. Now
$\langle \xi \rangle_{\vec{\ell}'}^s = \langle \xi \rangle_{\vec{\ell}'}^sm(\xi)^{-1}m(\xi)$,
and the Fourier multipliers
\begin{align*}
{\mathcal F}^{-1}_{\xi \to x}\bigl(\langle \xi \rangle_{\vec{\ell}'}^sm(\xi)^{-1}\cdot\textup{Id}_X\bigr){\mathcal F}_{x \to \xi} &: L_p(\R^d,X) \to L_p(\R^d,X), \\
{\mathcal F}^{-1}_{\xi \to x}m(\xi)\cdot\textup{Id}_X{\mathcal F}_{x \to \xi} &: W_p^{s;\vec{\ell}'}(\R^d,X) \to L_p(\R^d,X)
\end{align*}
are continuous. This shows that
$$
{\mathcal F}^{-1}_{\xi \to x}\langle \xi \rangle_{\vec{\ell}'}^s\cdot\textup{Id}_X{\mathcal F}_{x \to \xi} :
W_p^{s;\vec{\ell}'}(\R^d,X) \to L_p(\R^d,X)
$$
is continuous, thus giving the other inclusion $W_p^{s;\vec{\ell}'}(\R^d,X) \subset H_p^{s;\vec{\ell}'}(\R^d,X)$.
\end{remark}

\medskip

\noindent
In the case $\dim X < \infty$ and $\dim Y < \infty$, Theorem \ref{IterationRboundedness} below
strengthens classical results about norm estimates for pseudodifferential operators depending
on parameters (see \cite{Shubin}) to the $\RR$-boundedness of these families. It is
that theorem combined with a parametrix construction in the calculus of pseudodifferential
operators that yields the $\RR$-boundedness of resolvents $\lambda(A-\lambda)^{-1}$ of
anisotropic elliptic operators $A$ in Section \ref{Quasi}, and by a localization argument also
of elliptic scattering operators in Section \ref{Scattering}.

\begin{theorem}\label{IterationRboundedness}
Let $X$ and $Y$ be Banach spaces of class $\HT$ with property $(\alpha)$, and let
$a(x,\xi,\lambda) \in S^0_{\cl}(\R_x^d,S_{\cl}^{\mu;\vec{\ell}}(\R_{\xi}^d\times\Lambda;X,Y))$.
Write $\vec{\ell} = (\vec{\ell}',\vec{\ell}'') \in \N^{d+q}$, and let $\nu \geq \mu$ be fixed.

The family of pseudodifferential operators
$\op_x(a)(\lambda) : \S(\R^d,X) \to \S(\R^d,Y)$ extends by continuity to
$$
\op_x(a)(\lambda) : H_p^{s;\vec{\ell}'}(\R^d,X) \to H_p^{s-\nu;\vec{\ell}'}(\R^d,Y)
$$
for every $s \in \R$ and $1 < p < \infty$, and the operator function
$$
\Lambda \ni \lambda \mapsto \op_x(a)(\lambda) \in \L\bigl(H_p^{s;\vec{\ell}'}(\R^d,X),H_p^{s-\nu;\vec{\ell}'}(\R^d,Y)\bigr)
$$
belongs to the $\RR$-bounded symbol space $S_{\RR}^{\mu';\vec{\ell}''}(\Lambda;H_p^{s;\vec{\ell}'}(\R^d,X),H_p^{s-\nu;\vec{\ell}'}(\R^d,Y))$
with $\mu' = \mu$ if $\nu \geq 0$, or $\mu' = \mu-\nu$ if $\nu < 0$.

The mapping $\op_x : a(x,\xi,\lambda) \mapsto \op_x(a)(\lambda)$ is continuous in the symbol spaces
$$
S^0_{\cl}(\R_x^d,S_{\cl}^{\mu;\vec{\ell}}(\R_{\xi}^d\times\Lambda;X,Y)) \to
S_{\RR}^{\mu';\vec{\ell}''}(\Lambda;H_p^{s;\vec{\ell}'}(\R^d,X),H_p^{s-\nu;\vec{\ell}'}(\R^d,Y)).
$$
\end{theorem}
\begin{proof}
Let us begin with proving that
\begin{equation}\label{OpSobSpaces}
\op_x(a)(\lambda) : H_p^{s;\vec{\ell}'}(\R^d,X) \to H_p^{s-\nu;\vec{\ell}'}(\R^d,Y)
\end{equation}
is continuous, and that the set $\{\langle \lambda \rangle_{\vec{\ell}''}^{-\mu'}\op_x(a)(\lambda) \st \lambda \in \Lambda\}$
is an $\RR$-bounded subset of $\L\bigl(H_p^{s;\vec{\ell'}}(\R^d,X),H_p^{s-\nu;\vec{\ell'}}(\R^d,Y)\bigr)$
with $\RR$-bound dominated by $C\cdot p(a)$ with a constant $C \geq 0$ not depending on $a(x,\xi,\lambda)$
and a continuous seminorm $p$ on
$S^0_{\cl}(\R_x^d,S_{\cl}^{\mu;\vec{\ell}}(\R_{\xi}^d\times\Lambda;X,Y))$.

To this end note that
$$
S^0_{\cl}(\R_x^d,S_{\cl}^{\mu;\vec{\ell}}(\R_{\xi}^d\times\Lambda;X,Y)) \cong
S^0_{\cl}(\R_x^d)\hat{\otimes}_{\pi}S_{\cl}^{\mu;\vec{\ell}}(\R_{\xi}^d\times\Lambda;X,Y),
$$
hence it suffices to show this assertion for $\lambda$-dependent families of Fourier multipliers with symbols
$a(\xi,\lambda) \in S_{\cl}^{\mu;\vec{\ell}}(\R_{\xi}^d\times\Lambda;X,Y)$, and to
prove that multipliers
\begin{equation}\label{multipliers}
M(x)\textup{Id}_Y : H_p^{s-\nu;\vec{\ell}'}(\R^d,Y) \to H_p^{s-\nu;\vec{\ell}'}(\R^d,Y),
\quad M(x) \in S^0_{\cl}(\R_x^d),
\end{equation}
are continuous with continuous dependence on $M(x)$.

Let us consider the case of families of Fourier multipliers with symbols $a(\xi,\lambda)$,
and assume first that $\mu = \nu = 0$.
In view of $\langle D_x \rangle_{\vec{\ell}'}^{s}\op_x(a)(\lambda)\langle D_x \rangle_{\vec{\ell}'}^{-s}
= \op_x(a)(\lambda)$ the desired assertion for arbitrary $s \in \R$ reduces to $s = 0$. Let
$$
{\mathcal T} = \{\xi^{\beta}\partial_{\xi}^{\beta}a(\xi,\lambda)\st (\xi,\lambda) \in \R^d\times\Lambda,\,
\beta \leq (1,\ldots,1)\} \subset \L(X,Y).
$$
The mapping
$$
S_{\cl}^{0;\vec{\ell}}(\R_{\xi}^d\times\Lambda;X,Y) \to
S_{\cl}^{0;\vec{\ell}}(\R_{\xi}^d\times\Lambda;X,Y), \quad
a(\xi,\lambda) \mapsto \xi^{\beta}\partial_{\xi}^{\beta}a(\xi,\lambda)
$$
is continuous, and thus the set ${\mathcal T}$
is $\RR$-bounded in $\L(X,Y)$ by Proposition \ref{classicalRbounded} by $\tilde{C}\cdot q(a)$
with some constant $\tilde{C} \geq 0$ not depending on $a(\xi,\lambda)$ and a continuous seminorm
$q$ on $S_{\cl}^{0;\vec{\ell}}(\R_{\xi}^d\times\Lambda;X,Y)$.
Theorem \ref{iteratedFouriermultiplier} now gives the continuity of the Fourier multipliers
$\op_x(a)(\lambda) : L_p(\R^d,X) \to L_p(\R^d,Y)$ and the $\RR$-boundedness of the set
$$
\{\op_x(a)(\lambda)\st \lambda \in \Lambda\} \subset \L(L_p(\R^d,X),L_p(\R^d,Y))
$$
by a multiple of $q(a)$. This shows the assertion in the case $\mu = \nu = 0$.

Now consider the case of a general Fourier multiplier with symbol
$a(\xi,\lambda) \in S_{\cl}^{\mu;\vec{\ell}}(\R_{\xi}^d\times\Lambda;X,Y)$. Writing
$a(\xi,\lambda) = \langle \xi,\lambda \rangle_{\vec{\ell}}^{\mu}\cdot\bigl(\langle \xi,\lambda \rangle_{\vec{\ell}}^{-\mu}a(\xi,\lambda)\bigr)$
and noting that $a(\xi,\lambda) \mapsto \langle \xi,\lambda \rangle_{\vec{\ell}}^{-\mu}a(\xi,\lambda)$
is continuous in
$$
S_{\cl}^{\mu;\vec{\ell}}(\R_{\xi}^d\times\Lambda;X,Y) \to S_{\cl}^{0;\vec{\ell}}(\R_{\xi}^d\times\Lambda;X,Y)
$$
we obtain from the above that it is sufficient to show the $\RR$-boundedness of the family
$$
\langle \lambda \rangle_{\vec{\ell}''}^{-\mu'}\langle D_x,\lambda \rangle_{\vec{\ell}}^{\mu} :
H_p^{s;\vec{\ell}'}(\R^d,Y) \to H_p^{s-\nu;\vec{\ell}'}(\R^d,Y),
$$
i.e. the $\RR$-boundedness of the $\lambda$-dependent family of Fourier multipliers with symbols
$\langle \lambda \rangle_{\vec{\ell}''}^{-\mu'}\langle \xi,\lambda \rangle_{\vec{\ell}}^{\mu}\cdot\textup{Id}_Y$.
Evidently, this further reduces to consider the family of Fourier multipliers with symbols
$$
\langle \xi \rangle_{\vec{\ell}'}^{s-\nu}
\langle \lambda \rangle_{\vec{\ell}''}^{-\mu'}\langle \xi,\lambda \rangle_{\vec{\ell}}^{\mu}\langle \xi \rangle_{\vec{\ell}'}^{-s}\cdot\textup{Id}_Y
=
\langle \xi \rangle_{\vec{\ell}'}^{-\nu}
\langle \lambda \rangle_{\vec{\ell}''}^{-\mu'}\langle \xi,\lambda \rangle_{\vec{\ell}}^{\mu}\cdot\textup{Id}_Y
$$
on $L_p(\R^d,Y)$. By Kahane's contraction principle \eqref{Kahanecontraction} and
Theorem \ref{iteratedFouriermultiplier} it is sufficient to prove that the function
$\psi(\xi,\lambda) = \langle \xi \rangle_{\vec{\ell}'}^{-\nu}
\langle \lambda \rangle_{\vec{\ell}''}^{-\mu'}\langle \xi,\lambda \rangle_{\vec{\ell}}^{\mu}$
satisfies
$$
\sup\{|\xi^{\beta}\partial_{\xi}^{\beta}\psi(\xi,\lambda)|\st
(\xi,\lambda) \in \R^d\times\Lambda,\,\beta \leq (1,\ldots,1)\} < \infty.
$$
This, however, is an elementary estimate and follows easily.

We still have to consider the case of multipliers \eqref{multipliers}. Noting that
$$
S^{0}_{\cl}(\R_x^d) \ni M(x) \mapsto \langle \xi \rangle_{\vec{\ell}'}^{s-\nu}{\#}M(x){\#}\langle \xi \rangle_{\vec{\ell}'}^{-(s-\nu)}
\in S^{0}_{\cl}(\R_x^d,S_{\cl}^{0;\vec{\ell}'}(\R^d_{\xi}))
$$
is continuous, it suffices to show the continuity of pseudodifferential operators $\op_x\bigl(b(x,\xi)\cdot\textup{Id}_Y\bigr)$
on $L_p(\R^d,Y)$, where $b(x,\xi) \in S^{0}_{\cl}(\R_x^d,S_{\cl}^{0;\vec{\ell}'}(\R^d_{\xi}))$,
and the operator norm of $\op_x\bigl(b(x,\xi)\cdot\textup{Id}_Y\bigr)$ has to be bounded by a
multiple of a continuous seminorm of $b(x,\xi)$. As before write
$S^{0}_{\cl}(\R_x^d,S_{\cl}^{0;\vec{\ell}'}(\R^d_{\xi})) \cong
S^{0}_{\cl}(\R_x^d) \hat{\otimes}_{\pi} S_{\cl}^{0;\vec{\ell}'}(\R^d_{\xi})$,
which reduces the assertion to Fourier multipliers with symbols in
$S_{\cl}^{0;\vec{\ell}'}(\R^d_{\xi})\cdot\textup{Id}_Y$ and multipliers with symbols in
$S^{0}_{\cl}(\R_x^d)\cdot\textup{Id}_Y$ on $L_p(\R^d,Y)$. The case of Fourier multipliers
follows from Theorem \ref{iteratedFouriermultiplier}, and the case of multipliers on
$L_p(\R^d,Y)$ is elementary.

It remains to show that \eqref{OpSobSpaces} is indeed an $\RR$-bounded symbol of
(anisotropic) order $\mu' = \mu$ if $\nu \geq 0$ or
$\mu' = \mu-\nu$ if $\nu < 0$, respectively, with symbol estimates dominated by a continuous seminorm of
$a(x,\xi,\lambda) \in S^0_{\cl}(\R_x^d,S_{\cl}^{\mu;\vec{\ell}}(\R_{\xi}^d\times\Lambda;X,Y))$.
Note first that the mapping
$$
\Lambda \ni \lambda \mapsto a(x,\xi,\lambda) \in
S^0_{\cl}(\R_x^d,S_{\RR}^{\mu;\vec{\ell}'}(\R_{\xi}^d;X,Y))
$$
is $C^{\infty}$ (with $\lambda$-derivatives being represented by those of the symbol $a$), and
thus by what we have just proved we conclude that \eqref{OpSobSpaces} depends smoothly
on $\lambda \in \Lambda$, and
$$
\partial_{\lambda}^{\beta}\op_x(a)(\lambda) = \op_x\bigl(\partial_{\lambda}^{\beta}a\bigr)(\lambda) :
H_p^{s;\vec{\ell}'}(\R^d,X) \to H_p^{s-\nu;\vec{\ell}'}(\R^d,Y)
$$
for $\beta \in \N_0^q$. Applying the above we obtain that
$\{\langle \lambda \rangle_{\vec{\ell}''}^{-\mu'+|\beta|_{\vec{\ell}''}}\partial_{\lambda}^{\beta}\op_x(a)(\lambda)\st
\lambda \in \Lambda\}$ is $\RR$-bounded in $\L\bigl(H_p^{s;\vec{\ell}'}(\R^d,X),H_p^{s-\nu;\vec{\ell}'}(\R^d,Y)\bigr)$
by a multiple of a continuous seminorm of $\partial_{\xi}^{\beta}a(x,\xi,\lambda) \in 
S^0_{\cl}(\R_x^d,S_{\cl}^{\mu-|\beta|_{\vec{\ell}''};\vec{\ell}}(\R_{\xi}^d\times\Lambda;X,Y))$.
Since
$$
\partial_{\lambda}^{\beta} : S^0_{\cl}(\R_x^d,S_{\cl}^{\mu;\vec{\ell}}(\R_{\xi}^d\times\Lambda;X,Y)) \to
S^0_{\cl}(\R_x^d,S_{\cl}^{\mu-|\beta|_{\vec{\ell}''};\vec{\ell}}(\R_{\xi}^d\times\Lambda;X,Y))
$$
is continuous the proof of the theorem is complete.
\end{proof}

\section{Maximal regularity for anisotropic elliptic operators on $\R^d$}\label{Quasi}

\noindent
To illustrate how the results of the previous sections can be applied to prove maximal
regularity results for partial differential operators, we consider in this section the
case of anisotropic elliptic operators in $\R^d$.

Throughout this section let $X$ be a Banach space of class $\HT$ satisfying property $(\alpha)$.
Let $\vec{\ell}' \in \N^d$ be a vector determining the anisotropy, and let
\begin{equation}\label{QuasiPDO}
A = \sum\limits_{|\alpha|_{\vec{\ell}'} \leq \mu}a_{\alpha}(x)D_x^{\alpha} :
\S(\R^d,X) \to \S(\R^d,X),
\end{equation}
where the $a_{\alpha}(x) \in S^0_{\cl}(\R^d,\L(X))$ are operator valued coefficient functions, and $\mu \in \N$.

Let $\Lambda \subset \C$ be a closed sector. We assume the following anisotropic ellipticity
condition of $A$ with respect to $\Lambda$:

\begin{definition}\label{quasielliptparam}
$A$ is called \emph{parameter-dependent anisotropic elliptic with respect to $\Lambda$} if the
following two conditions are fulfilled:
\begin{enumerate}[i)]
\item The spectrum of the \emph{(anisotropic) principal symbol}
$$
\sum\limits_{|\alpha|_{\vec{\ell}'} = \mu}a_{\alpha}(x)\xi^{\alpha} \in \L(X)
$$
intersected with $\Lambda$ is empty for all $\xi \in \R^d\setminus\{0\}$ and all $x \in \R^d$.
\item The spectrum of the \emph{extended principal symbol}
$$
\sum\limits_{|\alpha|_{\vec{\ell}'} = \mu}a_{\alpha,(0)}(x)\xi^{\alpha} \in \L(X)
$$
intersected with $\Lambda$ is empty for all $\xi \in \R^d\setminus\{0\}$ and all
$x \in \R^d\setminus\{0\}$, where $a_{\alpha,(0)}(x)$ is the principal component of the
$\L(X)$-valued classical symbol $a_{\alpha}(x) \in S^0_{\cl}(\R^d,\L(X))$.
\end{enumerate}
The extended principal symbol (on $|x| = 1$) can be regarded as an extension of the
anisotropic principal symbol to the radial compactification of $\R^d$ (in the $x$-variables).
\end{definition}

\noindent
By the remarks given in the introduction, maximal regularity for anisotropic elliptic
operators $A$ (up to a spectral shift) follows from Theorem \ref{RbddResQuasi} below provided
that $A$ is parameter-dependent anisotropic elliptic with respect to the right half-plane
$\Lambda = \{\lambda \in \C \st \Re(\lambda) \geq 0\} \subset \C$, see Corollary \ref{MaxRegCorAniso}.
Note that the Sobolev spaces $H_p^{s;\vec{\ell'}}(\R^d,X)$ are of class $\HT$ and satisfy
property $(\alpha)$ in view of the permanence properties of these conditions
since they are isomorphic to $L_p(\R^d,X)$.

\begin{theorem}\label{RbddResQuasi}
Let $A$ be parameter-dependent anisotropic elliptic with respect to the closed sector $\Lambda \subset \C$.
Then
\begin{equation}\label{ASobSpCont}
A : H_p^{s+\mu;\vec{\ell}'}(\R^d,X) \to H_p^{s;\vec{\ell}'}(\R^d,X)
\end{equation}
is continuous for every $s \in \R$ and $1 < p < \infty$, and $A$ with domain
$H_p^{s+\mu;\vec{\ell}'}(\R^d,X)$ is a closed operator in $H_p^{s;\vec{\ell}'}(\R^d,X)$.

For $\lambda \in \Lambda$ with $|\lambda| \geq R$ sufficiently large the operator
$$
A - \lambda : H_p^{s+\mu;\vec{\ell}'}(\R^d,X) \to H_p^{s;\vec{\ell}'}(\R^d,X)
$$
is invertible for all $s \in \R$, and the resolvent
$$
\{\lambda(A-\lambda)^{-1}\st \lambda \in \Lambda,\,|\lambda| \geq R\} \subset
\L\bigl(H_p^{s;\vec{\ell}'}(\R^d,X)\bigr)
$$
is $\RR$-bounded.
\end{theorem}
\begin{proof}
The continuity of \eqref{ASobSpCont} follows from Theorem \ref{IterationRboundedness}.

The operator $A - \lambda : \S(\R^d,X) \to \S(\R^d,X)$ is of the form $\op_x(a)(\lambda)$
with the symbol
$$
a(x,\xi,\lambda) = \sum\limits_{|\alpha|_{\vec{\ell}'} \leq \mu}a_{\alpha}(x)\xi^{\alpha} - \lambda
\in S^{0}_{\cl}(\R^d_x,S^{\mu;\vec{\ell}}_{\cl}(\R^d_{\xi}\times\Lambda;X,X)),
$$
where $\vec{\ell} = (\vec{\ell}',\mu,\mu) \in \N^{d+2}$. Note that $\Lambda \subset \C \cong \R^2$
is considered as a real $2$-dimensional parameter space.

By our assumption of parameter-dependent anisotropic ellipticity, the parameter-de\-pen\-dent principal symbol
$$
a_{(\mu)}(x,\xi,\lambda) = \sum\limits_{|\alpha|_{\vec{\ell}'} = \mu}a_{\alpha}(x)\xi^{\alpha} - \lambda \in \L(X)
$$
is invertible for all $x \in \R^d$ and $(\xi,\lambda) \in \bigl(\R^d\times\Lambda\bigr)\setminus\{0\}$,
and with any excision function $\chi \in C^{\infty}(\R^{d+2})$ of the origin (i.e. $\chi \equiv 0$
near the origin and $\chi \equiv 1$ near infinity) we have
$$
b(x,\xi,\lambda) = \chi(\xi,\lambda)a_{(\mu)}(x,\xi,\lambda)^{-1}
\in S^{0}_{\cl}(\R^d_x,S^{-\mu;\vec{\ell}}_{\cl}(\R^d_{\xi}\times\Lambda;X,X)).
$$
We conclude that
$$
a{\#}b - 1,\; b{\#}a - 1 \in S^{0}_{\cl}(\R^d_x,S^{-1;\vec{\ell}}_{\cl}(\R^d_{\xi}\times\Lambda;X,X)),
$$
and the standard formal Neumann series argument now implies the existence of
\begin{align*}
p(x,\xi,\lambda) &\in S^{0}_{\cl}(\R^d_x,S^{-\mu;\vec{\ell}}_{\cl}(\R^d_{\xi}\times\Lambda;X,X)), \\
r_j(x,\xi,\lambda) &\in S^{0}_{\cl}(\R^d_x,S^{-\infty}(\R^d_{\xi}\times\Lambda;X,X)), \quad j=1,2,
\end{align*}
such that $a{\#}p = 1 + r_1$ and $p{\#}a = 1 + r_2$.

Let $P(\lambda) = \op_x(p)(\lambda)$, and $R_j(\lambda) = \op_x(r_j)(\lambda)$, $j=1,2$.
By Theorem \ref{IterationRboundedness} we have
\begin{align*}
P(\lambda) &\in S^{0;(\mu,\mu)}_{\RR}\bigl(\Lambda;H_p^{s;\vec{\ell}'},H_p^{s+\mu;\vec{\ell}'}\bigr) \cap
S^{-\mu;(\mu,\mu)}_{\RR}\bigl(\Lambda;H_p^{s;\vec{\ell}'},H_p^{s;\vec{\ell}'}\bigr), \\
R_j(\lambda) &\in \S\bigl(\Lambda,\L\bigl(H_p^{s;\vec{\ell}'},H_p^{t;\vec{\ell}'}\bigr)\bigr), \quad j=1,2,
\end{align*}
for all $s,t \in \R$. For $\lambda \in \Lambda$ with $|\lambda| \geq R$, the operators
$$
1 + R_j(\lambda) : H_p^{t;\vec{\ell}'}(\R^d,X) \to H_p^{t;\vec{\ell}'}(\R^d,X)
$$
are invertible for every $t \in \R$, and the inverses are represented as $1 + R_j'(\lambda)$
with
$$
R_j'(\lambda) \in \S\bigl(\Lambda,\L\bigl(H_p^{s;\vec{\ell}'}(\R^d,X),H_p^{s';\vec{\ell}'}(\R^d,X)\bigr)\bigr), \quad j=1,2,
$$
for any $s,s' \in \R$. More precisely, we may write
$$
R_j'(\lambda) = - R_j(\lambda) + R_j(\lambda)\chi(\lambda)\bigl(1+R_j(\lambda)\bigr)^{-1}R_j(\lambda)
$$
with some excision function $\chi \in C^{\infty}(\R^2)$ of the origin, and
$\bigl(1+R_j(\lambda)\bigr)^{-1}$ is the inverse of $1+R_j(\lambda)$ on some space
$H_p^{t_0;\vec{\ell}'}(\R^d,X)$ with a fixed $t_0 \in \R$.

We obtain that for $\lambda \in \Lambda$ with $|\lambda| \geq R$ the operator
$$
A - \lambda : H_p^{s+\mu;\vec{\ell}'}(\R^d,X) \to H_p^{s;\vec{\ell}'}(\R^d,X)
$$
is invertible for every $s \in \R$, and the resolvent is represented as
$$
(A-\lambda)^{-1} = P(\lambda) + P(\lambda)R_1'(\lambda).
$$
In view of Theorem \ref{IterationRboundedness} and Corollary \ref{BilderRbounded}
the proof is therefore complete.
\end{proof}

\begin{corollary}\label{MaxRegCorAniso}
Let $1 < p < \infty$, and assume that $A$ is parameter-dependent anisotropic elliptic
with respect to $\Lambda = \{\lambda \in \C \st \Re(\lambda) \geq 0\}$.
Then there exists $\gamma \in \R$ such that $A + \gamma$ with domain
$H_p^{s+\mu;\vec{\ell}'}(\R^d,X) \subset H_p^{s;\vec{\ell}'}(\R^d,X)$ has maximal
regularity for every $s \in \R$.
\end{corollary}

\section{Maximal regularity for elliptic scattering operators}\label{Scattering}

\noindent
Let $\Mbar$ be a $d$-dimensional smooth compact manifold with boundary. The aim of
this section is to prove maximal regularity for elliptic scattering operators on
$\Mbar$. These are elliptic differential operators in the interior $M = \open{\Mbar}$
which degenerate at the boundary in a specific way. The model example in
this context is $\Mbar = {\mathbb B}$, a $d$-dimensional ball. In this case an elliptic
scattering operator on ${\mathbb B}$ corresponds to an elliptic operator on $\R^d$
whose coefficients behave in some nice way when $|x| \to \infty$ radially, the ball
${\mathbb B}$ appears as the radial compactification of $\R^d$.

We follow Melrose \cite{Melrose} with our presentation of scattering
operators. The presence of coefficients in the morphisms of a bundle of (infinite dimensional)
Banach spaces represents no major difficulty as far as it concerns the action of the operators
on smooth sections that vanish to infinite order at the boundary of $\Mbar$.

\subsection*{Scattering differential operators}

Let $\x$ be a defining function for the boundary of $\Mbar$, i.e. $\x \in C^{\infty}(\Mbar)$
with $\x > 0$ on $M$, $\x = 0$ on $\partial\Mbar$, and $d\x \neq 0$ on $\partial\Mbar$.
Let ${}^{\textup{b}}\V(\Mbar)$ denote the vector fields on $\Mbar$ which are tangent on $\partial\Mbar$,
and let
$$
{}^{\sc}\V(\Mbar) = \x{}^{\textup{b}}\V(\Mbar)
$$
be the Lie algebra of \emph{scattering vector fields} on $\Mbar$. In coordinates
near $\partial\Mbar$, the vector fields in ${}^{\sc}\V(\Mbar)$ are spanned by
\begin{equation}\label{LocalscatteringVF}
x^2\frac{\partial}{\partial x}, \; x\frac{\partial}{\partial y_j}, \quad j = 1,\ldots,d-1,
\end{equation}
where $y_1,\ldots,y_{d-1}$ are coordinates on $\partial\Mbar$ and $x$ is a local boundary
defining function.

Let $\scT(\Mbar) \to \Mbar$ denote the \emph{scattering tangent bundle}, i.e. the vector bundle on
$\Mbar$ whose sections are the scattering vector fields. Fibrewise, we may represent
$\scT(\Mbar)$ as
$$
\scT_p(\Mbar) = {}^{\sc}\V(\Mbar)/{\mathcal I}_p(\Mbar)\cdot{}^{\sc}\V(\Mbar),
$$
where ${\mathcal I}_p(\Mbar) \subset C^{\infty}(\Mbar)$ is the ideal of functions that vanish
at $p$. Locally near $\partial\Mbar$, the vector fields \eqref{LocalscatteringVF} form a
smooth basis for $\scT(\Mbar)$. Let $\scT^*(\Mbar) \to \Mbar$ be the \emph{scattering cotangent bundle},
the dual of $\scT(\Mbar)$.

Let $\scDiff^*(\Mbar)$ be the enveloping algebra generated by ${}^{\sc}\V(\Mbar)$
and $C^{\infty}(\Mbar)$ consisting of the \emph{scattering differential operators}.
The operators of order $\mu \in \N_0$ are denoted as usual by $\scDiff^{\mu}(\Mbar)$.
The principal symbol $\sym(A)$ on $T^*M\setminus 0$ of an operator $A \in \scDiff^{\mu}(\Mbar)$
lifts to a well defined homogeneous function $\scsym(A)$ of degree $\mu$ on
$\scT^*\Mbar\setminus 0$ which is called the \emph{principal scattering}\ or
\emph{principal $\sc$-symbol} of $A$.
Note that $\scsym(A)$ is different from the symbol map considered in Section 6 of \cite{Melrose},
we are not heading for a Fredholm theory here.
Instead, $\scsym(A)$ may somewhat be regarded as a unified version of the principal symbol and
the extended principal symbol from Definition \ref{quasielliptparam} (when $\Mbar = {\mathbb B}$
is a ball).

\bigskip

\noindent
We are interested in maximal regularity of elliptic scattering operators with
operator valued coefficients. To this end, let $X \to \Mbar$ be a smooth vector bundle
of Banach spaces that (fibrewise) is of class $\HT$ and satisfies property $(\alpha)$
($X$ is the restriction of a corresponding bundle of Banach spaces with these properties from a
neighboring smooth manifold without boundary to $\Mbar$, e.g. from the double
$2\Mbar$). Note that both $\HT$ and property $(\alpha)$ are topological properties of
a Banach space, hence these notions are well defined for bundles of Banach spaces.
By considering the connected components of $\Mbar$ separately if necessary,
we may assume without loss of generality that $X$ can be locally trivialized with respect
to a fixed Banach space $X_0$, i.e. $X|_U \cong U \times X_0$ locally.

By $\L(X) \to \Mbar$ we denote the bundle of continuous linear operators in the fibres
of $X$. A scattering differential operator $A \in \scDiff^{\mu}(\Mbar;X)$ with coefficients
in $\L(X)$ is an operator of the form
\begin{equation}\label{scatteringoperatorvalued}
A = \sum\limits_{j=1}^N \varphi_j \otimes B_j : \dot{C}^{\infty}(\Mbar,X) \to \dot{C}^{\infty}(\Mbar,X),
\end{equation}
where $\varphi_j \in C^{\infty}(\Mbar,\L(X))$ and $B_j \in \scDiff^{\mu}(\Mbar)$, $j = 1,\ldots,N$, $N \in \N$,
and $\dot{C}^{\infty}(\Mbar,X)$ is the ($C^{\infty}(\Mbar)$-module) of smooth sections
of $X$ that vanish to infinite order on the boundary $\partial\Mbar$. The operator \eqref{scatteringoperatorvalued}
has an (evidently defined) principal $\sc$-symbol
$$
\scsym(A) \in C^{\infty}\bigl(\scT^*\Mbar \setminus 0, \L(\scpi^*X)\bigr),
$$
where $\scpi : \scT^*\Mbar \setminus 0 \to \Mbar$ is the canonical projection.

\subsection*{Function spaces and pseudodifferential operators}

Before coming to a general manifold $\Mbar$, we consider first the special case
$$
{\mathbb S}^d_+ = \{z'=(z'_1,\ldots,z'_{d+1}) \in \R^{d+1};\; |z'| = 1, \; z'_1 \geq 0\}.
$$
Consider the stereographic projection
$$
\SP : \R^d \ni z \longmapsto \Bigl(\frac{1}{(1+|z|^2)^{1/2}},\frac{z}{(1+|z|^2)^{1/2}}\Bigr) \in {\mathbb S}^d_+.
$$
For a Banach space $X_0$ of class $\HT$ having property $(\alpha)$ and $s \in \R$, $1 < p < \infty$,
the vector valued $\sc$-Sobolev space on ${\mathbb S}^d_+$ is defined as
$$
\scH_p^s({\mathbb S}^d_+,X_0) = \SP_* H_p^s(\R^d,X_0).
$$
Moreover, for a closed sector $\Lambda \subset \C$ and any $\ell \in \N$, we define the
class $\scPsi^{\mu;\ell}({\mathbb S}^d_+;\Lambda)$ as to consist of families of
operators
$$
A(\lambda) : \dot{C}^{\infty}({\mathbb S}^d_+,X_0) \to \dot{C}^{\infty}({\mathbb S}^d_+,X_0), \quad \lambda \in \Lambda,
$$
such that
$$
\bigl(\SP^*A(\lambda)\bigr)u(z) = (2\pi)^{-d}\int\limits_{\R^d}e^{iz\zeta}a(z,\zeta,\lambda)\hat{u}(\zeta)\,d\zeta, \quad u \in \S(\R^d,X_0),
$$
is a pseudodifferential operator with symbol
$$
a(z,\zeta,\lambda) \in S^0_{\cl}(\R_z^d,S_{\cl}^{\mu;\vec{\ell}}(\R_{\zeta}^d\times\Lambda;X_0,X_0)),
$$
where the vector $\vec{\ell}$ that determines the anisotropy of covariables $\zeta$ and parameters $\lambda$
is given by $\vec{\ell} = (\underbrace{1,\ldots,1}_{d},\ell,\ell) \in \N^{d+2}$. Note that $\Lambda \subset \C \cong \R^2$
is regarded as a real $2$-dimensional parameter space.

\bigskip

\noindent
In the general case of a compact manifold with boundary $\Mbar$ and a bundle $X \to \Mbar$
of Banach spaces of class $\HT$ having property $(\alpha)$, we define $\scH_p^s(\Mbar,X)$ as
to consist of all $u \in {\mathcal D}'(M,X)$ such that $\chi_*\bigl(\varphi u\bigr) \in
\scH_p^s({\mathbb S}^d_+,X_0)$ for all local charts $\chi : U \to {\mathbb S}^d_+$ with
$X|_U \cong U \times X_0$, and all $\varphi \in C^{\infty}(\Mbar)$ with compact support
contained in $U$.
The invariance of scattering pseudodifferential operators (see also below) and Theorem
\ref{IterationRboundedness} (a version of that theorem without parameters is
sufficient) imply that the spaces $\scH_p^s(\Mbar,X)$ are well defined for every $s \in \R$
and $1 < p < \infty$. Moreover, the projective topology with respect to the mappings
$$
\scH_p^s(\Mbar,X) \ni u \mapsto \chi_*\bigl(\varphi u\bigr) \in
\scH_p^s({\mathbb S}^d_+,X_0)
$$
for all charts $\chi$ and cut off functions $\varphi$ (as well as trivializations of $X$) makes
$\scH_p^s(\Mbar,X)$ a topological vector space which is normable so as to be a Banach space
that contains $\dot{C}^{\infty}(\Mbar,X)$ as dense subspace.

\medskip

\noindent
The class $\scPsi^{\mu;\ell}(\Mbar;\Lambda)$ of (anisotropic) parameter-dependent scattering
pseudodifferential operators on $\Mbar$ with coefficients in $\L(X)$ consists of operator
families
\begin{equation}\label{ScatteringCore}
A(\lambda) : \dot{C}^{\infty}(\Mbar,X) \to \dot{C}^{\infty}(\Mbar,X), \quad \lambda \in \Lambda,
\end{equation}
such that the following holds:
\begin{itemize}
\item For all $\varphi,\psi \in C^{\infty}(\Mbar)$ with disjoint supports we have
$$
\bigl(\varphi A(\lambda) \psi\bigr) u(z) = \int\limits_{M}k(z,z')u(z')\,\m(z'), \quad u \in \dot{C}^{\infty}(\Mbar,X),
$$
with $k(z,z') \in \dot{C}^{\infty}\bigl(\Mbar\times\Mbar,\L(\pi_L^*X,\pi_R^*X)\bigr)$, where
$\pi_L,\pi_R : \Mbar \times \Mbar \to \Mbar$ are the canonical projections on the left
and right factor, respectively, and $\m$ is any scattering density, i.e.
$\x^{d+1}\m$ is a smooth everywhere positive density on $\Mbar$ (recall that $d = \dim \Mbar$).
\item For any chart $\chi : U \to {\mathbb S}^d_+$, $U \subset \Mbar$, with
$X|_U \cong U \times X_0$, and all $\varphi,\psi \in C^{\infty}(\Mbar)$ with compact supports
contained in $U$, the operator push-forward $\chi_*\bigl(\varphi A(\lambda) \psi\bigr)$
is required to belong to the class $\scPsi^{\mu;\ell}({\mathbb S}^d_+;\Lambda)$ as defined
above.
\end{itemize}
A localization argument and Theorem \ref{IterationRboundedness} now imply the following

\begin{theorem}\label{IterationScattering}
Let $A(\lambda) \in \scPsi^{\mu;\ell}(\Mbar;\Lambda)$. Then \eqref{ScatteringCore} extends
by continuity to a family of continuous operators 
$$
A(\lambda) : \scH_p^s(\Mbar,X) \to \scH_p^{s-\nu}(\Mbar,X)
$$
for every $s \in \R$ and all $1 < p < \infty$, $\nu \geq \mu$. The operator function
$$
\Lambda \ni \lambda \mapsto A(\lambda) \in \L\bigl(\scH_p^s(\Mbar,X),\scH_p^{s-\nu}(\Mbar,X)\bigr)
$$
belongs to the $\RR$-bounded symbol space $S^{\mu';(\ell,\ell)}_{\RR}\bigl(\Lambda;\scH_p^s(\Mbar,X),\scH_p^{s-\nu}(\Mbar,X)\bigr)$
with $\mu' = \mu$ if $\nu \geq 0$, or $\mu' = \mu - \nu$ if $\nu < 0$.
\end{theorem}

\subsection*{$\RR$-boundedness of resolvents}

We are now ready to prove the $\RR$-boundedness of resolvents of elliptic scattering
differential operators with operator valued coefficients. As mentioned in the
introduction, maximal regularity (up to a spectral shift) is a consequence if we
choose $\Lambda \subset \C$ to be the right half-plane in Theorem~\ref{MaxRegScattering} below,
see Corollary \ref{MaxRegCorrScatt}.
Note that the $\sc$-Sobolev spaces are of class $\HT$ and satisfy property $(\alpha)$
as they are isomorphic to a finite direct sum of $X_0$-valued $L_p$-spaces, and $X_0$
has these properties.

\begin{theorem}\label{MaxRegScattering}
Let $\Lambda \subset \C$ be a closed sector, and let $A \in \scDiff^{\mu}(\Mbar;X)$, $\mu > 0$,
be a scattering differential operator on $\Mbar$ with coefficients in $\L(X)$, where
$X \to \Mbar$ is a smooth vector bundle of Banach spaces of class $\HT$ satisfying
property $(\alpha)$.
Assume that
$$
\spec\bigl(\scsym(A)(z,\zeta)\bigr) \cap \Lambda = \emptyset
$$
for all $(z,\zeta) \in \scT^*\Mbar \setminus 0$.

Then $A : \scH_p^{s+\mu}(\Mbar,X) \to \scH_p^{s}(\Mbar,X)$ is continuous for every $s \in \R$ and
$1 < p < \infty$, and $A$ with domain $\scH_p^{s+\mu}(\Mbar,X)$ is a closed operator
in $\scH_p^s(\Mbar,X)$.

Moreover, for $\lambda \in \Lambda$ with $|\lambda| \geq R$ sufficiently large, the operator
\begin{equation}\label{AmlScatSob}
A - \lambda : \scH_p^{s+\mu}(\Mbar,X) \to \scH_p^{s}(\Mbar,X)
\end{equation}
is invertible for all $s \in \R$, and the resolvent
$$
\{\lambda(A-\lambda)^{-1}\st \lambda \in \Lambda, \; |\lambda| \geq R\} \subset \L\bigl(\scH_p^s(\Mbar,X)\bigr)
$$
is $\RR$-bounded.
\end{theorem}
\begin{proof}
In view of the parameter-dependent ellipticity condition on the principal $\sc$-symbol of $A$
we conclude that we can construct local parametrices as in the proof of Theorem \ref{RbddResQuasi}
by symbolic inversion and a formal Neumann series argument. Patching these parametrices together
on $\Mbar$ with a partition of unity gives a global parameter-dependent
parametrix $P(\lambda) \in \scPsi^{-\mu;\mu}(\Mbar;\Lambda)$ of $A - \lambda$, i.e.
$$
\bigl(A-\lambda\bigr)P(\lambda) - 1, \; \bigl(A-\lambda\bigr)P(\lambda) - 1 \in
\scPsi^{-\infty;\mu}(\Mbar;\Lambda).
$$
By Theorem \ref{IterationScattering} we hence conclude that \eqref{AmlScatSob} is invertible
for all $s \in \R$ and $|\lambda| \geq R$ sufficiently large, and for these $\lambda$ 
we may write
$$
(A-\lambda)^{-1} - P(\lambda) \in \S\bigl(\Lambda,\L(\scH_p^{s}(\Mbar,X),\scH_p^{t}(\Mbar,X))\bigr)
$$
for all $s,t \in \R$. Theorem \ref{IterationScattering} and Corollary \ref{BilderRbounded} now
imply the assertion.
\end{proof}

\begin{corollary}\label{MaxRegCorrScatt}
Let $1 < p < \infty$, and let $A \in \scDiff^{\mu}(\Mbar;X)$, $\mu > 0$. We assume that
$$
\spec\bigl(\scsym(A)(z,\zeta)\bigr) \cap \{\lambda \in \C \st \Re(\lambda) \geq 0\} = \emptyset
$$
for all $(z,\zeta) \in \scT^*\Mbar \setminus 0$.
Then there exists $\gamma \in \R$ such that $A + \gamma$ with domain $\scH_p^{s+\mu}(\Mbar,X)
\subset \scH_p^{s}(\Mbar,X)$ has maximal regularity for every $s \in \R$.
\end{corollary}

\enlargethispage*{1cm}

\end{document}